\documentclass[11pt,reqno]{amsart}
\usepackage{graphicx, amssymb, hyperref, color,pdfsync}
\usepackage[all,cmtip]{xy}
\numberwithin{equation}{section}
\usepackage{mathrsfs}
\usepackage[symbol]{footmisc}

\usepackage{tikz}
\usetikzlibrary{matrix,arrows}
\usepackage{dsfont}
\usepackage[USenglish]{babel}
\usepackage{longtable}
\usepackage{tikz-cd}
\input epsf
\advance\hoffset by -0.9 truecm

\begin{document}
\newcommand{\s}{\vspace{0.2cm}}
\renewcommand{\thefootnote}{\alph{footnote}}

\newtheorem{theo}{Theorem} 
\newtheorem{prop}[theo]{Proposition}
\newtheorem{coro}[theo]{Corollary}
\newtheorem{lemm}[theo]{Lemma}
\newtheorem{rema}{Remark}

\theoremstyle{remark}
\newtheorem{defi}[theo]{\bf Definition}
\newtheorem{example}[theo]{\bf Example}

\title[On non-normal subvarieties of the moduli space]{On non-normal subvarieties of the moduli \\space of Riemann surfaces}
\date{}
\author[Hidalgo]{Rub\'en A. Hidalgo}
\address{Departamento de Matem\'atica y Estad\'istica, Universidad de La Frontera, Avenida Francisco Salazar 01145, Temuco, Chile.}
\email{ruben.hidalgo@ufrontera.cl}

\author[Paulhus]{Jennifer Paulhus}
\address{Department of Mathematics, Grinnell College, Grinnell, Iowa, U.S.A.}
\email{jenpaulhus@gmail.com}

\author[Reyes-Carocca]{Sebasti\'an Reyes-Carocca} 
\address{Departamento de Matem\'aticas, Facultad de Ciencias, Universidad de Chile, Las Palmeras 3425, \~{N}u\~noa, Santiago, Chile.}
\email{sebastianreyes.c@uchile.cl}
\author[Rojas]{Anita M. Rojas} 
\address{Departamento de Matem\'aticas, Facultad de Ciencias, Universidad de Chile, Las Palmeras 3425, \~{N}u\~noa, Santiago, Chile.}
\email{anirojas@uchile.cl}
\thanks{Corresponding author: Sebasti\'an Reyes-Carocca (sebastianreyes.c@uchile.cl)}

\begin{abstract}
In this article, we consider certain irreducible subvarieties of the moduli space of compact Riemann surfaces determined by the specification of actions of finite groups. We address the general problem of determining which among them are non-normal subvarieties of the moduli space. We obtain several new examples of subvarieties with this property. 
\end{abstract}

\maketitle
\thispagestyle{empty}
\section{Introduction}

Let $T_g$ denote the Teichm\"{u}ller space of compact Riemann surfaces of genus $g \geqslant 2.$ It is well-known that $T_g$ has the structure of a complex analytic manifold of dimension $3g-3$ and that the mapping class group $\mbox{Mod}_g$ acts properly discontinuously on $T_g$ by biholomorphisms. We denote by $$\Pi:T_g \to T_g/\mbox{Mod}_g$$the associated canonical projection. Since the group of (conformal) automorphisms of a compact Riemann surface of genus at least two is finite, a theorem due to Cartan (see \cite{Cartan}) implies that the moduli space $$\mathscr{M}_g:=T_g/\mbox{Mod}_g$$has the structure of a normal complex analytic space (complex orbifold) of dimension $3g-3.$ 

\s

If $t \in T_g$ is represented by the compact Riemann surface $S_t$, then the stabiliser of $t$ in $\mbox{Mod}_g$ can be identified with the group of (conformal) automorphisms  of $S_t,$ namely, $$\mbox{Stab}_{\tiny \mbox{Mod}_g}(t) \cong \mbox{Aut}(S_t).$$ 

For $g \geqslant 4$  the (orbifold) singular locus of $\mathscr{M}_g$ agrees with the branch locus of $\Pi$ 
$$\mbox{Sing}(\mathscr{M}_g)=\{[S] \in \mathscr{M}_g : \mbox{Aut}(S) \neq \{1\}\}.$$We refer to \cite{Popp} and \cite{Rauch}, and also  \cite{Nag} for more details.

\s

By contrast, the singular locus of $\mathscr{M}_2$ consists only of one point, algebraically represented by the Bolza curve $$y^2=x^5-1,$$whereas the singular locus of $\mathscr{M}_3$ consists of those points representing isomorphism classes of Riemann surfaces with non-trivial automorphisms, with the exception of those hyperelliptic whose reduced automorphism group is trivial (see \cite{Oort13}).

\s

In this work, we consider certain irreducible subvarieties of $\mathscr{M}_g$ that are determined by the specification of conjugacy classes of finite groups in the mapping class group. These subvarieties, which will be introduced in detail in \S\ref{preli}, can be understood in terms of Fuchsian groups as follows.
We denote by \begin{equation}\label{eq}\mathscr{M}_g(H, s, \theta) \end{equation}  the locus formed by the points of $\mathscr{M}_g$ representing isomorphism classes of compact Riemann surfaces of genus $g$ endowed with the action of a group $H$  and with topological class determined by  the surface-kernel epimorphism $$\theta : \Delta \to H,$$where $\Delta$ is a Fuchsian group of  signature $s$ (see \S\ref{preli} for more details). These loci have the structure of  irreducible subvarieties of $\mathscr{M}_g$ and $$\mbox{dim}(\mathscr{M}_g(H, s, \theta))=3h-3+r,$$  provided that $s=(h; m_1, \ldots, m_r).$ In general, such subvarieties are non-smooth and are contained in the singular locus of $\mathscr{M}_g$ (see \cite{B90} and \cite{GG92}).

\s

In this article, we deal with the general problem of deciding which subvarieties of $\mathscr{M}_g$ of type \eqref{eq} are non-normal. In addition, in the case of non-normality we address the problem of determining the set of non-normal points of the subvariety. 

\s

The concept of {\it normal variety} was introduced by Zariski in \cite{Z39} in the context of both affine and projective spaces. Following Mumford \cite{Mumford}, normality can be understood as a way to separate the ``branches'' of an algebraic variety at a singular point.

The importance of normal varieties lie in part in the fact that the co-dimension of their singular loci are greater than or equal to two (for instance, normal complex curves are smooth and the singularities of normal complex surfaces are isolated). Roughly speaking, a normal complex space only closely misses being a smooth complex manifold.

\s

The key fact we will employ to prove our results (and which will be explained in a precise way later) is that, as shown by Gonz\'alez-Diez and Harvey in \cite{GG92}, the existence of a non-normal point $$[X_0] \in \mathscr{M}_g(H, s, \theta)$$is intimately related to the existence of groups of automorphisms of $X_0$ isomorphic to $H$ acting on $X_0$ with topological class determined by $\theta$ and that are non-conjugate in the full automorphism group of $X_0$. 
\s

As we shall discuss later, for these subvarieties the normality can be understood as the property of being  biholomorphically equivalent to quotients of certain complex analytic submanifolds of $T_g$ by the action of appropriate groups of biholomorphisms of them.

\s

We mention some known facts concerning this problem, where the cyclic group of order $n$ is denoted by $C_n$, and the direct product $C_n \times \stackrel{m}{\ldots} \times C_n$ is denoted by $C_n^m$.

\s

{\bf 1.} The uniqueness of the hyperelliptic involution together with the results proved by Gonz\'alez-Diez in \cite{GG95}  show that if  $p$ is a prime number then $$\mathscr{M}_g(C_p, s=(0; p, \stackrel{r}{\ldots}, p), \theta)$$ is normal, for each action $\theta$. Here $ g=(r-2)(p-1)/2$ and $r \geqslant 3.$

\s

{\bf 2.} The first example of a non-normal subvariety of type \eqref{eq} associated to  regular covers of the projective line $\mathbb{P}^1$ was constructed by  Gonz\'alez-Diez and Hidalgo in  \cite{GG97}, being the two-dimensional subvariety of $\mathscr{M}_9$ given by $$\mathscr{M}_9(C_8, s=(0; 4,4,4,8,8), \theta),$$for some action $\theta$ suitably chosen.  This two-dimensional subvariety has a one-dimensional sublocus of non-normal points. Later, this construction was generalised by Carvacho in \cite{C13} to two-dimensional subvarieties associated to cyclic $2$-groups of the form $$\mathscr{M}_{3(2^n-1)}(C_{2^{n+1}}, s=(0; 2^n,2^n,2^n,2^{n+1},2^{n+1}), \theta_n) \, \mbox{ where } \, n \geqslant 2.$$

{\bf 3.} An example of a non-normal subvariety of type \eqref{eq} associated to  regular non-cyclic covers of $\mathbb{P}^1$ was found by Cirre in \cite{Ci04}, being the three-dimensional subvariety $$\mathscr{M}_3(C_2^2, s=(0; 2,2,2,2,2,2), \theta_h),$$where $\theta_h$ is the unique surface-kernel epimorphism for which $\mathscr{M}_3(C_2^2, s, \theta_h)$ lies inside the hyperelliptic locus of $\mathscr{M}_3$.

\s

{\bf 4.} For each integer $n \geqslant 3$ and for each prime number $p$ large enough, Hidalgo in \cite{Hidalgo11} gave a $((n-1)p-3)$-dimensional non-normal subvariety of the form $$\mathscr{M}_g(C_p^{n-1}, s=(0; p, \stackrel{(n-1)p}{\ldots}, p), \theta_{p,n}),$$for a suitable action $\theta_{p,n}.$ Here $g=1+p^{n-1}[(p-1)(n-1)-2]/2.$

\s

We should note that it is not difficult to construct examples of non-normal subvarieties of type \eqref{eq} with associated orbit space of positive genus. By contrast, to the best of our knowledge, apart from the aforementioned cases together with an example in genus three with $H \cong C_4$ (see \cite[Remark (1)]{GG97}) and an example in genus nine with $H \cong C_8$ (see \cite[Remark 4.3]{BRT23}), there are no more known examples of non-normal subvarieties of the desired type and associated to branched regular covers of $\mathbb{P}^1.$ 

\s

This paper is aimed at providing new examples of subvarieties of this type. The results, which will be stated in \S\ref{resultados}, can be briefly summarised as follows.

\s

{\bf a.} We provide the full list of non-normal subvarieties in $\mathscr{M}_2$ and $\mathscr{M}_3$ of type \eqref{eq}. In addition, we describe explicitly the set of their non-normal points (Theorems \ref{caseg2} and \ref{caseg3}). 

\s

{\bf b.} We prove that, for each $g \geqslant 2,$ the $g$-dimensional subvariety  of $\mathscr{M}_g$$$\mathscr{M}_g(C_2^2, s=(0; 2, \stackrel{g+3}{\ldots},2), \theta_h),$$where $\theta_h$ stands for the action corresponding to the hyperelliptic locus, is non-normal. We also describe a set of non-normal points of it (Theorem \ref{funciona}). This generalises fact ${\bf 3}$ above.

\s

{\bf c.} We deduce that, for each $g\geqslant 2,$ the moduli space $\mathscr{M}_g$ contains a non-normal subvariety of type \eqref{eq} where $s$ is a genus zero signature (Corollary \ref{todog}).

\s

{\bf d.} We prove that, for each $n \geqslant 1,$ the dihedral group of order $8n$ yields a non-normal subvariety of the moduli space $\mathscr{M}_{2n+1}$. The novelty of this result is that it produces the first non-abelian examples of type \eqref{eq} where $s$ is a genus zero signature (Theorem \ref{conj2}).

\s

{\bf e.} We prove that, for each odd integer $n \geqslant 1,$ the cyclic group of order $2n$ gives rise to a non-normal subvariety of the moduli space $\mathscr{M}_{(n-1)^2}$. This result can be seen as an extension of  fact  {\bf 2} above from cyclic $2$-groups to cyclic groups of even order (Theorem \ref{conj7}).

\s

{\bf f.} We show that,  for each pair of integers $n \geqslant 4$ and $k \geqslant 2$, the family of the so-called generalised Fermat curves gives rise to a non-normal subvariety of $\mathscr{M}_g$ of the form $$\mathscr{M}_g(C_k^{n-1}, s=(0; k,\stackrel{k(n-1)}{\ldots},k), \theta_{k,n}),$$for some suitably chosen action $\theta_{k,n}$. Here $g=1+k^{n-1}[(n-1)(k-1)-2]/2.$ This result (Theorem \ref{CFG}) extends fact ${\bf 4}$ above.

\section{Preliminaries} \label{preli}
\subsection*{Fuchsian groups and  Riemann surfaces} Let $\mathbb{H}$ be the upper half-plane and let $\Delta$ be a co-compact Fuchsian group. The  {\it signature} of $\Delta$ is the tuple \begin{equation} \label{sig} s=s(\Delta)=(h; m_1, \ldots, m_l),\end{equation}where $h$ denotes the genus of the quotient surface $\mathbb{H}/{\Delta}$ and $m_1, \ldots, m_l$ the branch indices in the associated universal projection $\mathbb{H} \to \mathbb{H}/{\Delta}.$ If $l=0$ then it is said that $\Delta$ is a {\it surface} Fuchsian group, and if $h=0$ and $l=3$ then $\Delta$ is  called a {\it triangle} Fuchsian group.

\s

If $\Delta$ is a Fuchsian group of signature \eqref{sig} then $\Delta$ has a canonical presentation
\begin{equation*}\label{prese}\langle \alpha_1, \ldots, \alpha_{h}, \beta_1, \ldots, \beta_{h}, x_1, \ldots , x_l : x_1^{m_1}, \cdots, x_l^{m_l}, \Pi_{i=1}^{h}[\alpha_i, \beta_i] \Pi_{j=1}^l x_j\rangle,\end{equation*}where the brackets stands for the commutator. 

\s

Let $S$ be a compact Riemann surface of genus $g \geqslant 2$. By the uniformisation theorem, there is a surface Fuchsian group $\Gamma$ such that $$S \cong \mathbb{H}/{\Gamma}.$$Furthermore, by Riemann's existence theorem,  a finite group $H$ acts on $S$ if and only if there is a Fuchsian group $\Delta$ together with a group  epimorphism \begin{equation*}\label{epi}\theta: \Delta \to H \, \, \mbox{ such that }  \, \, \mbox{ker}(\theta)=\Gamma.\end{equation*}In such a case, the genus $g$ of $S$ is related to $s(\Delta)$ by the Riemann-Hurwitz formula $$2(g-1)=|H|(2(h-1)+\Sigma_{j=1}^l(1-1/m_j)).$$

We say that $H$ acts on $S$ with signature $s(\Delta)$ and that this action is  represented by the {\it surface-kernel epimorphism} $\theta$  (from now on, simply SKE). We usually identify $\theta$ with the tuple or {\it generating vector} $$\theta=(\theta(\alpha_1), \ldots, \theta(\alpha_h), \theta(\beta_1), \ldots, \theta(\beta_h), \theta(x_1), \ldots, \theta(x_l)) \in H^{2h+l}.$$

Let $G$ be a group and let $H$ be a subgroup of $G.$ The action of $H$ on $S$ is said to {\it tightly extend} to an action of $G$ if the following three statements hold.\begin{enumerate}
\item There is a Fuchsian group $\Delta'$ with $\Delta \leqslant \Delta'$.
\item There is a surface-kernel  epimorphism $$\Theta: \Delta' \to G \, \, \mbox{ with }  \, \, \Theta|_{\Delta}=\theta.$$
\item The Teichm\"{u}ller spaces of $\Delta$ and $\Delta'$ have the same dimension.
\end{enumerate} 
An action is called {\it maximal} if it cannot be tightly extended. See \cite{Bu03}, \cite{Ries93} and \cite{Sing72}. We also refer to the survey article \cite{BPW22} for more details and related problems concerning group actions on Riemann surfaces.

\subsection*{The subvariety $\mathscr{M}_g(H, s, \theta)$}

Let $S$ be a  Riemann surface of genus $g \geqslant 2$ and let  $$\psi:H \to \mbox{Aut}(S)$$ be an action of a finite group $H$ on $S$, with signature $s.$ Assume that $S'$ is another Riemann surface of genus $g$ endowed with an action $\psi':H \to \mbox{Aut}(S')$. If there exists a (orientation-preserving) homeomorphism  \begin{equation} \label{lapiz12}\phi: S \to S' \mbox{ such that } \phi \psi(H) \phi^{-1}=\psi'(H)\end{equation}then we say that the actions of $H$  are in the same topological class or that they are {\it topologically equivalent}. Clearly, in such a case the  signatures of the actions agree. 

\s

The topological equivalence of actions can be understood in terms of Fuchsian groups as follows. If we write $$S/{H} \cong \mathbb{H}/{\Delta} \, \mbox{ and } \, S'/H \cong \mathbb{H}/{\Delta}'$$then each (orientation-preserving)  homeomorphism $\phi$ as in \eqref{lapiz12}  induces  a group isomorphism $\phi^*: \Delta \to \Delta'$. Hence, we may assume $\Delta=\Delta'$. We denote the subgroup of $\mbox{Aut}(\Delta)$ consisting of such automorphisms $\phi^*$ by $\mathscr{B}$. Following \cite{B90}, the SKEs $$\theta : \Delta \to H \, \mbox{ and } \, \theta': \Delta \to H$$representing two actions of $H$  are topologically equivalent if and only if there are $\omega \in \mbox{Aut}(H)$ and  $\phi^* \in \mathscr{B}$ such that  $$\theta' = \omega\circ\theta \circ \phi^*.$$ 

Sources for the characterisation of topological actions by certain purely algebraic data include Nielsen \cite{Nielsen}, Harvey \cite{H71} and Gilman \cite{Gilman}. 
\s

As mentioned in the introduction, we denote by $$\mathscr{M}_g(H, s, \theta) \subset \mathscr{M}_g$$ the locus in $\mathscr{M}_g$ consisting of the points representing isomorphism classes of Riemann surfaces endowed with an action topologically equivalent to the one of $H$ on $S.$ According to \cite{B90}, this locus is a closed irreducible (non necessarily smooth) subvariety of  $\mathscr{M}_g$. 

\s

We also recall that the topological class of the action of $H$ is  determined by a group monomorphism $\iota: H \to \mbox{Mod}_g$, and that $$\mathscr{M}_g(H, s, \theta) = \Pi(T_g(H, s, \theta)),$$where $T_g(H, s, \theta)$ stands for the fixed locus  of $\iota(H)$ in $T_g$.

\s

We refer to the articles \cite{BCI14} and \cite{CI10}  for the case of low genus, and to \cite{BRT23} for algorithms to find concrete examples of subvarieties of this kind. See also \cite{Pign}.

 \s
 
We remark that in some papers the notation $\bar{\mathscr{M}}_g(H, s, \theta)$ is used instead of ${\mathscr{M}}_g(H, s, \theta).$

\subsection*{The normalisation of $\mathscr{M}_g(H, s, \theta)$ and a criterion for non-normality}Let $S'$ and $S''$ be two Riemann surfaces of genus $g \geqslant 2$ endowed with actions $$\psi': H \to \mbox{Aut}(S') \, \mbox{ and } \, \psi'': H \to \mbox{Aut}(S'').$$

Assume that $\psi'$ and $\psi''$ are topologically equivalent to the action of $H$ on $S$ (of signature $s$ and represented by the SKE $\theta$). The actions are termed {\it analytically (or conformally) equivalent} if there exists a conformal isomorphism  \begin{equation*} \label{eccc}\Phi: S' \to S'' \mbox{ such that } \Phi \psi'(H) \Phi^{-1}=\psi''(H).\end{equation*}

Observe that two actions on a Riemann surface $S$ are analytically equivalent if and only if the corresponding groups of automorphisms are conjugate in the full automorphism group of the surface.

\s

We denote by $\{(S', H, \psi')\}$ the analytic equivalence class determined by $\psi'$, and by $$\tilde{\mathscr{M}}_g(H, s, \theta)$$the set of such classes. Observe that the natural correspondence $$\pi : \tilde{\mathscr{M}}_g(H, s, \theta) \to {\mathscr{M}}_g(H, s, \theta) \subset \mathscr{M}_g \,\mbox{ given by } \{(S', H, \psi')\} \mapsto [S']$$is surjective. 

Furthermore,  following the discussion in \cite[\S1]{GG92}, classical results of Teichm\"{u}ller theory imply that  $$\tilde{\mathscr{M}}_g(H, s, \theta) \cong T_g(H, s, \theta)/ \mbox{Stab}_{\tiny \mbox{Mod}_g}(T_g(H, s, \theta))$$and hence has the structure of a normal complex analytic space which makes $\pi$ a morphism of  complex analytic spaces. Indeed, the surjection $\pi$ turns out to be the normalisation of $\mathscr{M}_g(H, s, \theta)$. See \cite[Theorem 1]{GG92}. We have the following commutative diagram$$\begin{tikzcd}
  T_g(H,s, \theta) \arrow[r] \arrow[dr, "\Pi"]
    & \bar{\mathscr{M}}_g(H,s, \theta)\arrow[d, "\pi"]\\
&\mathscr{M}_g(H,s, \theta) \end{tikzcd}$$

Observe, in addition, that  $\pi$ fails to be a biholomorphism if and only if  $\mathscr{M}_g(H, s, \theta)$ is non-normal. In other words, the following statements are equivalent.

\s

{\bf 1.}  $\mathscr{M}_g(H, s, \theta)$ is a non-normal subvariety of $\mathscr{M}_g.$

\s

{\bf 2.} There is $t \in T_{g}(H, s, \theta)$ and  $\varphi \in \mbox{Mod}_g$ such that $\varphi(t) \in T_{g}(H, s, \theta)$ but $$\varphi \notin\mbox{Stab}_{\tiny \mbox{Mod}_g}(T_g(H, s, \theta)).$$In such a case both $H$ and $\varphi H \varphi^{-1}$ fix $t$ but they are different.

\s

{\bf 3.} There is  $[X_0] \in \mathscr{M}_g(H, s, \theta)$ with an action of $H$ that is topologically but not analytically equivalent to the action of $H$ determined by $\theta$.

\s

We emphasise here that if the action $\theta$ of $H$ on Riemann surfaces of genus $g$ is maximal, then the implication $$S \in \mathscr{M}_g(H, s, \theta) \, \implies \, \mbox{Aut}(S) \cong H$$holds for each $S$ with the possible exception of some Riemann surfaces that necessarily lie in a positive codimensional sublocus of $\mathscr{M}_g(H, s, \theta)$. See \cite{B90} for more details.

\begin{rema}\label{clona}
It is worth remarking that possibly some non-normal points of $\mathscr{M}_g(H, s, \theta)$ do not satisfy the statements {\bf 2} and {\bf 3} above. This situation is explained by the fact that the set of non-normal points of $\mathscr{M}_g(H, s, \theta)$ form a closed subset. We refer to \cite[\S2]{GG92} for more details.
\end{rema}

\subsection*{Notations} Let $n \geqslant 2$ be an integer. A  primitive $n$-th root of unity is denoted by $\omega_n$, the cyclic group of order $n$ is  denoted by $C_n,$ the dihedral group of order $2n$ is  denoted by $\mathbf{D}_n$, and the symbol $a^n$ in a signature abbreviates the expression $a, \stackrel{n}{\ldots},a.$ The pair $(n,i)$ used as a {\it group id} for a group of order $n$ refers to its label in the {\it small group} database employed by several computer algebra programs (SageMath, Magma and GAP, among others).

\section{Statement of the results}\label{resultados}

\begin{theo}\label{caseg2}
There is exactly one non-normal irreducible subvariety of $\mathscr{M}_2$ of type $\mathscr{M}_2(H, s, \theta).$ This subvariety has dimension two and is formed by the Riemann surfaces with a group of automorphisms isomorphic to $$H=C_2^2 \mbox{ acting with signature }  s=(0; 2^5).$$In addition, its set of non-normal points is formed by the Riemann surfaces with a group of automorphisms isomorphic to $$\mathbf{D}_4 \mbox{ acting with signature }  (0; 2^3,4).$$ In particular, the Accola-Maclachlan curve and the Wiman curve of type II $$y^2=x^6-1 \,\, \mbox{ and }\,\, y^2=x(x^4-1)$$respectively, are non-normal points.
\end{theo}

If we denote by $\mathfrak{B}$ the Bolza curve (namely, the unique compact Riemann surface of genus two endowed with an automorphism of order five) then the theorem above says that$$\{\mbox{branch locus of }T_2 \to \mathscr{M}_2\}-\{\mathfrak{B}\}$$is a non-normal subvariety of $\mathscr{M}_2.$ We recall that except for $\mathfrak{B}$ (which is an isolated point), all Riemann surfaces in the branch locus of $\mathscr{M}_2$ admit the (unique) action of  $C_2^2$.

\begin{theo}\label{caseg3}
There are exactly seven non-normal irreducible subvarieties of  $\mathscr{M}_3$ of type $\mathscr{M}_3(H, s, \theta)$. 

\s

{\bf 1.} The one-dimensional subvariety of Riemann surfaces with a group of automorphisms isomorphic to $$H=\mathbf{S}_4 \mbox{ acting with signature }  s=(0; 2^3,3).$$The subvariety has only one non-normal point, being the Klein quartic. 

\s

{\bf 2.} The two-dimensional subvariety of Riemann surfaces with a group of automorphisms isomorphic to $$H=C_2^3 \mbox{ acting with signature }  s=(0; 2^5).$$The set of non-normal points is formed by the Riemann surfaces with a group of automorphisms isomorphic to $$C_2 \times \mathbf{D}_4 \mbox{ acting with signature }  (0; 2^3,4).$$In particular, the most symmetric among the hyperelliptic Riemann surfaces of genus three is a non-normal point.

\s

{\bf 3.} The two-dimensional subvariety of Riemann surfaces with a group of automorphisms isomorphic to $$H=\mathbf{D}_4 \mbox{ acting with signature }  s=(0; 2^5).$$The set of non-normal points is formed by the Riemann surfaces with a group of automorphisms isomorphic to $$(C_4 \times C_2)\rtimes C_2 \mbox{ acting with signature }  (0; 2^3,4).$$In particular, the Fermat quartic is a non-normal point.

\s

{\bf 4.} The two-dimensional subvariety of Riemann surfaces with a group of automorphisms isomorphic to $$H=C_4 \mbox{ acting with signature }  s=(0; 2^3, 4^2).$$The set of non-normal points is formed by the Riemann surfaces with a group of automorphisms isomorphic to $$C_4 \times C_2 \mbox{ acting with signature }  (0; 2^2,4^3).$$In particular, the Wiman curve of type II is a non-normal point.

\s

{\bf 5.} The two irreducible components of the three-dimensional family of Riemann surfaces with a group of automorphisms isomorphic to $$H= C_2^2=\langle a,b : a^2, b^2, (ab)^2\rangle \mbox{ acting with signature }  s=(0; 2^6).$$

Such components are represented by the SKE $$(a,a,a,a,b,b) \, \mbox{ and } \, (a,a,b,b,ab,ab)$$and their sets of non-normal points are formed by the Riemann surfaces with a group of automorphisms isomorphic to$$C_2^3 \, \mbox{ and }  \mathbf{D}_4  \mbox{ acting with signature }  (0; 2^5)$$respectively.

\s

{\bf 6.} The four-dimensional subvariety of Riemann surfaces with a group of automorphisms isomorphic to $$H=C_2 \mbox{ acting with signature }  s=(1; 2^4).$$The set of non-normal points is formed by the Riemann surfaces with a group of automorphisms isomorphic to $$C_2^2=\langle a,b: a^2,b^2,(ab)^2\rangle \mbox{ acting with signature }  (0; 2^6)$$and represented by the SKE $(a,a,b,b,ab,ab).$
\end{theo}

\begin{theo}\label{funciona}
For each integer $g \geqslant 3$, the $g$-dimensional  subvariety of $\mathscr{M}_g$\begin{equation}\label{rss}\mathscr{M}_g(C_2^2, s=(0;2^{g+3}), \theta_h)\end{equation}is non-normal, where $\theta_h$ stands for  the unique SKE for which $\mathscr{M}_g(C_2^2, s, \theta_h)$ lies inside the hyperelliptic locus of $\mathscr{M}_g$.  The subvariety \eqref{rss} is formed by the surfaces $S$ represented by the algebraic curves $$y^2=\Pi_{i=1}^{g+1}(x-a_i)(x-\tfrac{1}{a_i})$$where $a_1, \ldots, a_{g+1}$ are pairwise distinct complex numbers satisfying that $a_ia_j \neq 1$ for all $i,j \in \{1, \ldots, g+1\}$. The group of automorphisms of $S$ that is isomorphic to $C_2^2$ is represented by  $$\langle (x,y)\mapsto (x,-y), \, (x,y)\mapsto (1/x,y/x^{g+1}) \rangle.$$

Furthermore, a set of non-normal points for \eqref{rss} is given as follows.

\s

{\bf a.} If $g$ is odd then a set of non-normal points consists of the  Riemann surfaces with a group of automorphisms isomorphic to $$C_2 \times \mathbf{D}_{g+1} \, \mbox{ acting with signature }\, (0; 2^3, g+1),$$and represented by the algebraic curves$$y^2=(x^{g+1}-a)(x^{g+1}-\tfrac{1}{a}) \, \mbox{ where } a \neq 0, \pm 1.$$In particular, the Accola-Maclachlan curve is a non-normal point of \eqref{rss}.

\s

{\bf b.} If $g$ is even then a set of non-normal points consists of the Riemann  surfaces with a group of automorphisms isomorphic to $$\mathbf{D}_{g} \, \mbox{ acting with signature }\, (0; 2^4, g),$$and  represented by the algebraic curves$$y^2=x(x^{\frac{g}{2}}-a)(x^{\frac{g}{2}}-\tfrac{1}{a})(x^{\frac{g}{2}}-b)(x^{\frac{g}{2}}-\tfrac{1}{b})$$where  $a,b \neq 0, \pm 1$ are distinct and $ab \neq 1.$ In particular, the Wiman curve of type II is a non-normal point of \eqref{rss}.
\end{theo}

The next result follows directly from Theorems \ref{caseg2} and \ref{funciona}. 

\begin{coro}\label{todog} There are infinitely many compact Riemann surfaces  of genus $g \geqslant 2$ that are branched regular covers of $\mathbb{P}^1$  endowed with two analytically non-conjugate groups of automorphisms whose actions are topologically equivalent.
\end{coro}

\begin{theo}\label{conj2} 
For each odd integer $g \geqslant 3$, the two-dimensional subvariety $$\mathscr{M}_g(\mathbf{D}_{2(g-1)}, s=(0; 2^5), \theta)$$is non-normal, where $\theta$ stands for the SKE representing the unique action of $\mathbf{D}_{2(g-1)}$ in genus $g$ with signature $s.$ A set of non-normal points is formed by the Riemann surfaces with a group of automorphisms isomorphic to a semidirect product$$C_{2(g-1)} \rtimes C_2^2 \, \mbox{ acting with signature } (0; 2^3, 4).$$
\end{theo}

\begin{theo}\label{conj7} 
For each odd integer $n \geqslant 3$, the moduli space $\mathscr{M}_{(n-1)^2}$ contains a non-normal subvariety of dimension $n-1$ of the form $$\mathscr{M}_{(n-1)^2}(C_{2n}, s=(0; 2^2,n, \stackrel{n}{\ldots}, n), \theta).$$A set of non-normal points is formed by the Riemann surfaces with a group of automorphisms isomorphic to$$C_{2} \times C_{2n} \, \mbox{ acting with signature } (0; 2^2, n, \stackrel{(n-1)/2}{\ldots}, n,2n).$$
\end{theo}

\begin{theo}\label{CFG} Let $k \geqslant 2$ and $n \geqslant 4$ be integers. The subvariety $$\mathscr{M}_g(C_k^{n-1}, s=(0; k,\stackrel{k(n-1)}{\ldots},k), \theta_{k,n}) \, \mbox{ where } g=1+\tfrac{1}{2}k^{n-1}((n-1)(k-1)-2)$$is non-normal, for some action $\theta_{k,n}$. A set of non-normal points is formed by the so-called generalised Fermat curves of type $(k,n)$, algebraically represented by
$$ \left. \begin{array}{ccccccc}
x_{1}^{k}&+&x_{2}^{k}&+&x_{3}^{k}&=&0\\
\lambda_{1}x_{1}^{k}&+&x_{2}^{k}&+&x_{4}^{k}&=&0\\
\, &\vdots &\,&\vdots&\,&\vdots&\,\\
\lambda_{n-2}x_{1}^{k}&+&x_{2}^{k}&+&x_{n+1}^{k}&=&0
\end{array}\right \} \subset \mathbb{P}^n
$$ where $\lambda_1, \ldots, \lambda_{n-2}$ are pairwise distinct complex numbers different from 0 and 1. 
\end{theo}

\section{Proof of Theorem \ref{caseg2}}
The nontrivial groups of automorphisms of compact Riemann surfaces of genus two and the description of their actions  are well-known: there are exactly twenty possibilities (see, for instance \cite{Bre}, \cite{B91},  \cite{Conder}, \cite{K} and \cite{LMFDB}). However, among such  possibilities, only seven correspond to maximal actions. The following table summarises such  groups and their actions (the trivial case $\mbox{Aut}(S) \cong C_2$ has not been considered).

\s

\begin{center}
\begin{tabular}{|ccccc|} \hline
 label &  notation & signature & group & id \\ \hline \hline
 5 & $\mathscr{C}_5$& $ (0; 2^5)$ & $C_2^2$ & $(4,2)$ \\ \hline
 11 &$\mathscr{C}_{11}$& $ (0; 2^3, 4)$ & $\mathbf{D}_4$ & $(8,3)$ \\ \hline 
 14 &$X_{14}$& $ (0; 2, 5, 10)$ & $C_{10}$ & $(10,2)$ \\ \hline 
 15 &$\mathscr{C}_{15}$& $ (0; 2^3, 3)$ & $\mathbf{D}_{6}$ & $(12,4)$ \\ \hline 
 20 &$X_{20}$& $ (0; 2, 4, 6)$ & $(C_{6}\times C_2)\rtimes C_2$ & $(24,8)$ \\ \hline 
 21 &$X_{21}$& $ (0; 2, 3, 8)$ & $\mathbf{GL}_{2}(3)$ & $(48,29)$ \\ \hline 
\end{tabular}
\end{center}
$$\mbox{Table 1: Full automorphism groups of Riemann surfaces of genus $g=2$}$$

The enumeration of the cases is as given in \cite{K}.  We denote by $\mathscr{C}_i$ the family of Riemann surfaces corresponding to the $i$-th case. If $\mathscr{C}_i$ consists of only one point, then we  denote it simply by $X_i.$ We recall that for each possibility in the table above there is only one topological class of actions.  

\s

Assume that $S$ belongs to the one-dimensional family $\mathscr{C}_{15}$  (respectively, to $\mathscr{C}_{11}$), and that $\mathbf{D}_6$ (respectively,  that $\mathbf{D}_4$) is not isomorphic to the automorphism group of $S$. It follows that the order of the automorphism group of $S$ is a multiple of 12 and therefore $S$ is isomorphic to $X_{20}$ or $X_{21}.$ Observe that the automorphism group of each one of these Riemann surfaces contains, up to conjugation, only one subgroup isomorphic to $\mathbf{D}_6$ and only one subgroup isomorphic to $\mathbf{D}_4$. Hence,  none of the members of $\mathscr{C}_{15}$  (respectively,  $\mathscr{C}_{11}$) enjoys the property of having two non-conjugate subgroups isomorphic to $\mathbf{D}_6$ (respectively,  to $\mathbf{D}_4$). Therefore, the families $\mathscr{C}_{15}$ and $\mathscr{C}_{11}$ are normal subvarieties of $\mathscr{M}_2$.

\s

Observe that the intersection $\mathscr{C}_{15} \cap \mathscr{C}_{11}$ is formed by $X_{20}$ and $X_{21}$.

\s 

Assume now that $S$ belongs to the one-dimensional family $\mathscr{C}_{5}$ and that $C_2^2$ is not isomorphic to its automorphism group. By arguing as before we see that $S \in \mathscr{C}_{11} \cup \mathscr{C}_{15}.$ Again the fact that, up to conjugation, the groups $\mathbf{D}_6$ and $\mathbf{D}_4$ have one and two subgroups isomorphic to $C_2^2$ respectively implies that:

\begin{enumerate}
\item $\mathscr{C}_{5}$ is a non-normal subvariety of $\mathscr{M}_2$, and
\item the set of non-normal points consists of the family $\mathscr{C}_{11}$. 
\end{enumerate}Observe that, in particular, $X_{20}$ and $X_{21}$ are non-normal points of $\mathscr{C}_5$ as well.

\begin{rema}\mbox{}\label{telefono}
\begin{enumerate}
\item It is worth remarking that $X_{21}$ is a non-normal point of $\mathscr{C}_5$ in spite of the fact that 
its automorphism group has, up to conjugation, only one subgroup isomorphic to $C_2^2$. This follows from the fact that the set of non-normal points is closed, and that each $ X \in \mathscr{C}_{11}-\{X_{21}\}$ is a non-normal point of $\mathscr{C}_{5}$.   This is an explicit example of the situation anticipated in Remark \ref{clona}.

\s

\item An algebraic description for the members of $\mathscr{C}_{5}$ is given by the curves $$Y_{a,b}: y^2=(x^2-1)(x^2-a^2)(x^2-b^2)$$where $a \neq b$ are complex numbers different from $\pm 1$. The non-normal points of $\mathscr{C}_{5}$ are given by the curves $Y_{a,1/a}.$ 
\end{enumerate}
\end{rema}

\section{Proof of Theorem \ref{caseg3}}
The nontrivial groups of automorphisms of compact Riemann surfaces of genus three and their signatures  are well-known; see, for instance \cite{Bre}, \cite{B91},  \cite{K} and \cite{LMFDB}. Among such possibilities, only twenty-two correspond to full automorphism groups (see, for instance, \cite{S02}). These groups are summarised in the following table.

\s

\begin{center}
\begin{tabular}{|c c c c || c c c c |} 
 \hline
 case & sign & group &  id   & case & sign & group & id \\ [0.5ex] 
 \hline \hline
 3    & $(1; {2^4})$  & $C_2$ & $(2,1)$  &  26   & $(0; {2^3,6})$ &$\mathbf{D}_6$& $(12,4)$\\ \hline
 4    & $(0; {2^8})$   &  $C_2$& $(2,1)$  &   30   & $(0; {2, 7,14})$ & $C_{14}$& $(14,2)$\\ \hline
 6    & $(0; {3^5})$   &  $C_3$& $(3,1)$    &   31   &$(0; {2^3,4})$ & $C_2 \times \mathbf{D}_4$& $(16,11)$\\ \hline
9   & $(0; {2^6})$   & $C_2^2$&  $(4,2)$  &  32   & $(0; {2^3,4})$ & $(C_4 \times C_2)\rtimes C_2$& $(16,13)$\\ \hline
 10    & $(0; 2^3,4^2)$    &  $C_4$& $(4,1)$ & 38   & $(0; {2^3,3})$   & $\mathbf{S}_{4}$& $(24,12)$\\ \hline
 13    & $(0; {2^4,3})$   & $\mathbf{S}_3$& $(6,1)$ &   42  &  $(0; {2,4,12})$ &$C_{4} \times \mathbf{S}_3$& $(24,5)$\\ \hline
 14   & $(0; {2,3^2,6})$    & $C_6$& $(6,2)$&  43   & $(0; {2,4,8})$ & $(C_8 \times C_2)\rtimes C_2$& $(32,9)$\\ \hline
 19  & $(0; {2^5})$ &$\mathbf{D}_4$& $(8,3)$& 46 &  $(0; {2,4,6})$ & $C_2 \times \mathbf{S}_4$& $(48,48)$\\  \hline
 20   & $(0; {2^5})$   &$C_2^3$& $(8,5)$&  47   & $(0; {2, 3,12})$ & $\mbox{SL}(2,3) \rtimes C_2$& $(48,33)$\\ \hline
 21  &$(0; 2^2,4^2)$& $C_2 \times C_4$& $(8,2)$&  {48}   & $(0; {2,3,8})$ &$C_{4}^2 \rtimes \mathbf{S}_3$& $(96,64)$\\  \hline
 24   & $(0; 3,9^2)$  & $C_9$& $(9,1)$&  {49}   & $(0; {2,3,7})$  & $\mathbb{P}\mbox{SL}(3,2)$& $(168,42)$\\ \hline
\end{tabular}
$$\mbox{Table 2: Full automorphism groups of Riemann surfaces of genus $g=3$}$$
\end{center} As in the genus two case, the enumeration of the cases is as given in \cite{K}, we denote by $\mathscr{C}_i$ the family of Riemann surfaces corresponding to the $i$-th case, and if $\mathscr{C}_i$ consists of only one point then we denote it simply by $X_i.$

\s

The family $\mathscr{C}_4$ corresponds to the hyperelliptic case and so it is normal, whereas the family $\mathscr{C}_6$ is normal by the results of \cite{GG95}.

\s

We now proceed to study each positive dimensional case separately.

\s

{\bf a.} We consider the family $\mathscr{C}_{31}.$ If $S$ belongs to this family and $H=C_2 \times \mathbf{D}_4$ is not isomorphic to its automorphism group $G$  then the order of $G$ is a multiple of 16 and  acts on $S$ with a triangle signature. It follows that $S$ is isomorphic to one of $$X_{43}, \, X_{46}, \, X_{47} \mbox{ or } X_{48}.$$ 

The third and fourth cases do not need to be considered since the automorphism groups of these surfaces do not contain any subgroup isomorphic to $H$. In addition, the automorphism groups of $X_{43}$ and $X_{46}$ each contain only one conjugacy class of subgroups isomorphic to $H$. It then follows that there is no member of  $\mathscr{C}_{31}$ possessing two non-conjugate subgroups isomorphic to $H.$ Hence $\mathscr{C}_{31}$ is a normal subvariety of $\mathscr{M}_3$.

\s

By proceeding similarly, one sees that $\mathscr{C}_{14}, \mathscr{C}_{26}$ and $\mathscr{C}_{32}$ are normal subvarieties of $\mathscr{M}_3.$

\s

{\bf b.} We consider the family $\mathscr{C}_{21}.$ Assume that $S \in \mathscr{C}_{21}$ and write $$H=C_4 \times C_2=\langle a, b: a^4, b^2, aba^{-1}b\rangle$$ 

{\it Claim.} The action of $H$ on $S$ is represented by the SKE $(a^2,b,a,ab).$

\s

We recall that the group $H$ acts in genus three with signature $s=(0; 2^2, 4^2)$ in three different topological ways. Such actions are represented by the SKEs $$\theta_1=(a^2, b, a, ab), \, \theta_2=(b,b,a,a^{-1}) \, \mbox{ and } \, \theta_3=(b,a^2b,a,a).$$

We denote  the subvariety of Riemann surfaces $X$ determined by the SKE $\theta_i$ by $\mathcal{F}_i$ for $i=1,2,3$. Observe that $\mathcal{F}_1$ and $\mathcal{F}_2$ consist of hyperelliptic surfaces;  the hyperelliptic involution being represented by $a^2$ and $b$ respectively.

\s

Following \cite{Sing72}, if the action of $H$ tightly extends to a group $H'$ then necessarily $H'$ has order 16 and the  signature of the action of $H'$ is $s'=(0; 2^3,4).$  It then follows that either $X \in \mathscr{C}_{31}$ or $X\in\mathscr{C}_{32}.$ 

\s

Assume that $X \in  \mathscr{C}_{31}.$ Then $$H'=C_2 \times \mathbf{D}_4=\langle t, r, s: t^2,r^4,s^2,(sr)^2,trtr^{-1},(ts)^2\rangle$$and, since there is a unique topological action of this group with signature $s'$, we can assume that the action of $H'$ on $X$ is represented by the SKE $$\Theta=(s,srt, tr^2, r).$$ The group $H'$ has a unique subgroup isomorphic to $H$, namely,  $\tilde{H}=\langle r, t \rangle.$ Observe that the induced action of $\tilde{H}$ on $X$ is represented by the SKE $$\Theta|_{\tilde{H}}=(tr^2,tr^2,r,r^{-1}) \, \mbox{ which is equivalent to }\theta_2.$$We have then obtained that $\mathcal{F}_2=\mathscr{C}_{31}.$

\s

Assume that $X \in \mathscr{C}_{32}.$ Then $$H' \cong (C_4 \times C_2)\rtimes C_2.$$

There is a unique topological class of actions of $H'$ with signature $s'$. This uniqueness coupled with the fact that this group has subgroups isomorphic to $H$ show that  $$\mathscr{C}_{32}=\mathcal{F}_1 \, \mbox{ or } \, \mathscr{C}_{32}=\mathcal{F}_3.$$ 

The former case is impossible, since $X$ is non-hyperelliptic \cite{S02}. Hence $\mathscr{C}_{32}=\mathcal{F}_3.$

\s

The above arguments imply that the action represented by $\theta_1$ is maximal, and this proves the claim.

\s

We now assume that $H$ is not isomorphic to the automorphism group of $S \in \mathscr{C}_{21}=\mathcal{F}_1.$  Then the order of this group is a multiple of eight and it acts on $S$ with a triangle signature. Among the groups satisfying these properties, there are only two of them possessing, up to conjugation, more than one subgroup isomorphic to $H$. These groups give rise to the following possibilities for non-normal points of  $\mathscr{C}_{21}$: $$X_{43} \mbox{ (the Accola-Maclachlan curve) and } X_{48} \mbox{ (the Fermat quartic).}$$

The latter case must be disregarded since the Fermat quartic is non-hyperelliptic. On the other hand, the action of $$\mbox{Aut}(X_{48})\cong \langle x,y,z : x^8,y^2,z^2,[x,y], [z,y], zxzxy\rangle$$ on $X_{48}$ is represented by the SKE $\Theta=(z,zx,x^{-1}),$ and this group has exactly two subgroups $$H_1=\langle x^2,y \rangle \, \mbox{ and } \, H_2=\langle xz,x^4\rangle$$ isomorphic to $H.$ A computation shows that $$\Theta|_{H_1}=(y,y,x^2,x^6) \, \mbox{ and } \, \Theta|_{H_2}=(y,x^4,xz,x^5z).$$Therefore, the actions of $H_1$ and $H_2$ are topologically non-equivalent.

\s

We then conclude that $\mathscr{C}_{21}$ is a normal subvariety of $\mathscr{M}_3$. In a very similar way it can be seen that $\mathscr{C}_{13}$ is also a normal subvariety.

\s

{\bf c.} We consider the family $\mathscr{C}_{10}.$ A routine computation shows that among the automorphism group of Riemann surfaces of genus three, there are exactly two of them that have, up to conjugation, at least two subgroups isomorphic to $C_4$ in such a way that the induced actions have signature $s=(0;2^3,4^2).$ These groups give rise to the family $\mathscr{C}_{21}$ and  the surface $X_{42}$ (the Wiman curve of type II). 

\s

We claim that $X_{42}$ belongs to the family $\mathscr{C}_{21}.$ Indeed, the action of $$\mbox{Aut}(X_{42}) \cong \langle t,r,s:t^4,r^3,s^2,(sr)^2,stst^{-1}, rtr^{-1}t^{-1}\rangle$$on $X_{42}$ is represented by the SKE $$\Theta=(sr,st^{-1},tr).$$

The restricted action of the subgroup $\langle s, t \rangle \cong C_2 \times C_4$ is represented by the SKE $$\Theta|_{\langle s, t \rangle}=(s,st,t^2, t^{-1})$$which is equivalent to the SKE $\theta_1$ (introduced in item {\bf b}); this SKE represents $\mathscr{C}_{21}.$  

\s

The fact that there is only one topological class of actions of $C_4$ in genus three with signature $s$ allows us to assert that each member of $\mathscr{C}_{21}$ is endowed with two subgroups isomorphic to  $C_4$ whose actions are topologically but not analytically equivalent. Hence, we are in a position to conclude that $\mathscr{C}_{10}$ is a non-normal subvariety of $\mathscr{M}_3$ and that its set of non-normal points consists of the family $\mathscr{C}_{21}$. 

\s

By proceeding similarly, one can conclude that $\mathscr{C}_{19}, \mathscr{C}_{20}$ and $\mathscr{C}_{38}$ are non-normal subvarieties of $\mathscr{M}_3.$

\s

{\bf d.} We now consider the family $\mathscr{C}_9.$ This family consists of two irreducible components. In fact, if we write $$C_2^2=\langle a,b: a^2,b^2,(ab)^2\rangle$$then such components are $\mathscr{C}_9^1$ and $\mathscr{C}_9^2$ represented by $$\theta_1=(a,a,a,a,b,b) \, \mbox{ and } \,  \theta_2=(a,a,b,b,ab,ab)$$ respectively. We first observe that among the automorphism groups of Riemann surfaces of genus three there are exactly nine of them that have, up to conjugation, at least two subgroups isomorphic to $C_2^2$ acting with signature $s=(0;2^6).$ These groups yield the following possibilities for non-normal points of $\mathscr{C}_{9}$: \begin{equation}\label{chau}\mathscr{C}_{19}, \mathscr{C}_{20}, \mathscr{C}_{31}, \mathscr{C}_{32}, \mathscr{C}_{38}  \, X_{43}, X_{46}, {X}_{48}, X_{49}.\end{equation}

Assume that $S\in \mathscr{C}_{19}$. If we write $$H'=\langle r,s: r^4,s^2,(sr)^2\rangle \cong \mathbf{D}_{4}$$then the action of this group on $S$ is represented by the SKE $$\Theta=(sr,sr,s,sr^2,r^2).$$ This group has two conjugacy classes of subgroups isomorphic to $C_2^2,$ represented by $$H_1=\langle s,r^2\rangle \mbox{ and } H_2=\langle sr, r^2\rangle.$$

A computation shows that the action of $H_1$ and of $H_2$ is represented by the SKEs $$\Theta|_{H_1}=(s,s,sr^2,sr^2,r^2,r^2) \mbox{ and } \Theta|_{H_2}=(sr,sr,sr^3,sr^3,r^2,r^2),$$showing that  they are topologically equivalent to $\theta_2.$ We then conclude that $\mathscr{C}_{19}$  consists of non-normal points of $\mathscr{C}_9^2$.

\s

Assume that $S\in \mathscr{C}_{20}$. If we write $$H'=\langle x,y,z: x^2,y^2,z^2,(xy)^2,(xz)^2,(yz)^2\rangle \cong C_{2}^3$$then the action of this group on $S$ is represented by the SKE $$\Theta=(x,x,y,yz,z).$$ This group has four conjugacy classes of subgroups isomorphic to $C_2^2$ with induced action of signature $(0;2^6).$ Such classes are represented by $$H_1=\langle x,y\rangle, \, H_2=\langle x,z \rangle, \, H_3=\langle x,yz \rangle \mbox{ and } H_4=\langle y,z\rangle.$$ By proceeding as in the previous case, it can be seen that $$\Theta|_{H_1} \cong \Theta|_{H_2} \cong \Theta|_{H_3} \cong \theta_1 \mbox{ and } \Theta|_{H_4} \cong \theta_2.$$ Hence $\mathscr{C}_{20}$  consists of non-normal points of $\mathscr{C}_9^1$.

\s

For the remaining possibilities \eqref{chau} of non-normal points $S$ of $\mathscr{C}_{9}$, a routine computation shows that each $S$ has a group of automorphisms isomorphic to $\mathbf{D}_4$ or $C_2^2$ with signature $(0;2^5).$ Since these groups act in genus three with this signature in a unique topological way, we conclude that  $S$ belongs to either $\mathscr{C}_{19}$ or $\mathscr{C}_{20}$ and so we can ignore those remaining cases. Consequently, the set of non-normal of points of $\mathscr{C}_9^1$ agrees with $\mathscr{C}_{20}$ and the set of non-normal of points of $\mathscr{C}_9^2$ agrees with $\mathscr{C}_{19}.$

\s

{\bf e.} We finally consider the family $\mathscr{C}_3.$ The automorphism groups of Riemann surfaces of genus three that have, up to conjugation, at least two subgroups isomorphic to $C_2$ acting with signature $s=(1;2^4)$ give rise to the following possibilities for non-normal points: $$\mathscr{C}_9^2, \mathscr{C}_{19}, \mathscr{C}_{20}, \mathscr{C}_{31}, \mathscr{C}_{32}, \mathscr{C}_{38}\, X_{43}, X_{46}, {X}_{47}, X_{48}.$$

Observe that if $S \in \mathscr{C}_9^2$ then $S$ has three groups of automorphisms isomorphic to $C_2$ acting with signature $(1;2)$, showing that $\mathscr{C}_9^2$ consists of non-normal points of $\mathscr{C}_3.$ Each one of the remaining possibilities is contained in $\mathscr{C}_9^2$ and consequently the set of non-normal points of $\mathscr{C}_3$ agrees with $\mathscr{C}_9^2$.

\section{Proof of Theorem \ref{funciona}}
\subsection*{The case $g$ odd} Let $g \geqslant 3$ be an odd integer. The group $$G=C_2 \times \mathbf{D}_{g+1}=\langle t, r,s : t^2,r^{g+1},s^2,(sr)^2,(ts)^2,trtr^{-1} \rangle$$acts on a one-dimensional family of Riemann surfaces $\mathscr{N}_g$ with signature $(0; 2^3, g+1)$. In fact, the existence of the family follows from the SKE $$\theta=(t,tsr,s,r).$$The Riemann-Hurwitz formula says that the genus of such  surfaces is $g.$ Since the signature of the action is maximal, up to finitely many exceptions, the full automorphism group of the members of $\mathscr{N}_g$ is isomorphic to $G.$ Observe that the subgroups of $G$ $$H_1=\langle t,s \rangle \, \mbox{ and } \, H_2=\langle t,z \rangle \, \mbox{ where }z:=r^{\frac{g+1}{2}}$$are isomorphic but non-conjugate (the latter one is central). 

\s

We claim that the induced actions of $H_1$ and $H_2$ on each $S \in \mathscr{N}_g$ are topologically equivalent. To prove the claim we introduce the following notation: $$\pi : S \to S/G \, \mbox{ and }\, \pi_i : S \to S/H_i$$ denote the canonical projections given by the action of $G$ and $H_i$ on $S$ for $i=1,2$, and $q_1, q_2, q_3,q_4$ denote the ordered branch values of  $\pi$. 

\s

The following statements hold.

\s

{\bf a.} The $\pi$-fiber of $q_1$ consists of $2(g+1)$ branch points, all of them with $G$-stabiliser $\langle t \rangle.$ Since $\langle t \rangle \subset H_i$, these points have $H_i$-stabilisers $\langle t \rangle$ and give rise to $g+1$ branch values of $\pi_i$ for $i=1,2.$  

\s

{\bf b.}  The $\pi$-fiber of $q_2$ consists of $2(g+1)$ branch points with $G$-stabiliser of the form $\langle tsr^{m} \rangle$ where $m$ is odd. Since $\langle tsr^{m} \rangle$ and $H_i$ intersects trivially for each $m$, these points have $H_i$-stabilisers trivial and do not give rise to branch values of $\pi_i$ for $i=1,2.$  

\s

{\bf c.}  The $\pi$-fiber of $q_3$ consists of $2(g+1)$ branch points with $G$-stabiliser of the form $\langle sr^{n} \rangle$ where $n$ is even. Note that for each $n$, there are exactly 4 points with $G$-stabiliser  $\langle sr^{n} \rangle.$ Since $\langle sr^{n} \rangle \cap H_1 \neq \{1\}$ if and only if $n=0$ and $\langle sr^{n} \rangle \cap H_2 = \{1\}$ for all $n,$ one sees that exactly 4 of these points have non-trivial $H_1$-stabiliser  and all of them have trivial $H_2$-stabiliser. Thus, the $\pi$-fiber of $q_3$ gives rise to 2 branch values of $\pi_1$ with $H_1$-stabiliser $\langle s \rangle$, and does not give rise to branch values of $\pi_2.$ 

\s

{\bf d.}  The $\pi$-fiber of $q_4$ consists of $4$ branch points, all of them with $G$-stabiliser $\langle r \rangle.$  Since $\langle r \rangle \cap H_1 = \{1\}$ and $\langle r \rangle \cap H_2 = \langle z\rangle$, these points do not give rise to branch points of $\pi_1$ and give rise to 2 branch values of $\pi_2$, with $H_2$-stabiliser $\langle z\rangle$.

\s

All the above show that the signature of the induced action of $H_1$ and $H_2$ on $S$ is $(0;2^{g+3})$ each (see \cite[\S3]{Anita}). In addition,  if the induced action of $H_1$ is represented by $$\theta_1=(h_1, \ldots, h_{g+3})$$ then $\langle h_i \rangle$ must be $H_1$-conjugate to $\langle t \rangle$ for $1 \leqslant i \leqslant g+1$, and  $\langle h_i \rangle$ must be $H_1$-conjugate to $\langle s \rangle$ for $i=g+2, g+3.$ The fact that $H_1$ is abelian shows that necessarily$$\theta_1=(t, \stackrel{g+1}{\ldots},t, s,s).$$In an analogous way,  the induced action of $H_2$ on $S$ is represented by the SKE $$\theta_2=(t, \stackrel{g+1}{\ldots},t, z,z).$$Consequently, these actions are topologically equivalent, as claimed.

\s

Observe that the actions $\theta_1 \cong \theta_2$  correspond to the equisymmetric stratum $\theta_h$ of hyperelliptic Riemann surfaces with an action of $C_2^2$ with signature $(0;2^{g+3})$; the hyperelliptic involution being represented by $t$. Thus, we are in a position to conclude that \begin{equation}\label{fff}\mathscr{M}_g(C_2^2, s=(0; 2^{g+3}), \theta_h)\end{equation}contains the family $\mathscr{N}_g$ which enjoys the following property: up to finitely many exceptions, each $X \in \mathscr{N}_g$ has two isomorphic non-conjugate groups of automorphisms acting in a topologically equivalent way. In other words, \eqref{fff} is a non-normal subvariety of $\mathscr{M}_g$ and  $\mathscr{N}_g$ consists of non-normal points of it. 

\s

We now proceed to provide an algebraic description of \eqref{fff}. Let $S$ be a member of \eqref{fff} and let \begin{equation}\label{ppp}a_1, \ldots, a_{g+1}, b_1, \ldots, b_{g+1}\end{equation}denote the $2g+2$ branch values of the covering map $$S \to R=S/\langle t \rangle \cong \mathbb{P}^1$$induced by the hyperelliptic involution. Since $S$ has an involution $s$ that induces an involution on $\tilde{s}$ on $R$ keeping the set of points \eqref{ppp} invariant and the signature of the action of $H_1=\langle t,s \rangle$ on $S$ is $(0; 2^{g+3})$, we may assume  $$\tilde{s}(\xi)=\tfrac{1}{\xi} \, \mbox{ and therefore } b_i=\tfrac{1}{a_i} \mbox{ for } i \in \{1, \ldots, g+1\}.$$Thus, the surface $S$ is algebraically represented by the curve $$y^2=\Pi_{i=1}^{g+1}(x-a_i)(x-\tfrac{1}{a_i})$$where $a_ia_j \neq 1$ for each $i,j \in \{1, \ldots, g+1\}.$ In this model, the group $H_1$ is generated by the   transformations $$\langle t(x,y)=(x,-y), s(x,y)=(1/x,y/x^{g+1}) \rangle $$ 

Assume now that $S$ belongs to  $\mathscr{N}_g.$ Then $r$ induces an automorphism  $\tilde{r}$ of $R$ of order $g+1$ that remains the set of points \eqref{ppp} invariant. Since the signature of the action of $\langle t, r \rangle$ on $S$ is $(0; 2^2, g+1, g+1)$ one sees that, up to a M\"{o}bius transformation, $$\tilde{r}(\xi)=\omega_{g+1} \xi \, \mbox{ and therefore that } a_i= \omega_{g+1}^{i-1} a_1 \mbox { for } i=\{1, \ldots, g+1\}.$$Note that $\tilde{r}$ has $0$ and $\infty$ as fixed points.

Set $a:=a_1^{g+1}.$  The regular covering map induced by $\langle \tilde{r} \rangle$  $$R \to R'=R/\langle \tilde{r} \rangle \cong S/\langle t, r \rangle\,\mbox{ is given by } \,  \xi \mapsto \xi^{g+1},$$showing that the branch values of the composed map $S \to R'$ are $0, \infty$ (marked with $g+1$ each) and $a, 1/a$ (marked with $2$ each). 

\s

In turn, as $\langle t, r \rangle$ is a normal subgroup of $G$, the automorphism $s$ of $S$ induces an involution $\hat{s}$ of $R'$ in such a way that $S/G \cong R'/\langle \hat{s}\rangle.$ Clearly $\hat{s}(\xi)=1/\xi$ and $\pm 1$ are its fixed points. It then follows that the regular covering map induced by $\hat{s}$ $$R' \to R'/\langle \hat{s} \rangle \cong S/G \, \mbox{ is given by } \, \xi \mapsto \xi + 1/\xi,$$and the branch values of the composed map $S \to S/G$ are $\infty$ marked with $g+1$, $\pm 2$ marked with 2, and $a+1/a$ marked with 2.
 
\s

All the above says that each member of $\mathscr{N}_g$ is  represented by the  curve \begin{equation}\label{pz}y^2=\Pi_{i=1}^{g+1}(x-a_i)(x-\tfrac{1}{a_i})=(x^{g+1}-a)(x^{g+1}-\tfrac{1}{a})\end{equation}for some $a \neq 0, \pm 1.$ In this model, the group of automorphisms isomorphic to $G$ is generated by the transformations $$t(x,y)=(x,-y), \,s(x,y)=(1/x, y/x^{g+1}) \mbox{ and } r(x,y)=(\omega_{g+1} x, y).$$ 

\s

Finally, by taking $\lambda=\omega_4$ in \eqref{pz} one sees that the Accola-Maclachlan curve $$y^2=x^{2g+2}+1$$belongs to $\mathscr{N}_g,$ as desired.

\subsection*{The case $g$ even} Let $g \geqslant 4$ be an even integer. Consider the group $$G= \mathbf{D}_{g}=\langle r,s : r^{g},s^2,(sr)^2\rangle$$ and set $z=r^{\frac{g}{2}}$. The Riemann-Hurwitz formula and the SKE $$ \theta=(z,z, s, sr, r^{-1})$$guarantee the existence of a two-dimensional family $\mathscr{S}_g$ of Riemann surfaces of genus $g$ with a group of automorphisms isomorphic to $G$ acting with signature $(0; 2^4, g).$ Since the signature of the action  is maximal, generically, the full automorphism group of the members of $\mathscr{S}_g$ is isomorphic to $G.$ Observe that the subgroups of $G$ $$H_1=\langle s,z \rangle \, \mbox{ and } \, H_2=\langle  sr,z \rangle$$are isomorphic but non-conjugate. By proceeding analogously as in the previous case, one sees that the signature of the induced actions of $H_1$ and $H_2$ on $S$ is $(0;2^{g+3})$, and that they are represented by the SKE $$\theta_1=(z, \stackrel{g}{\ldots},z, s, sz, z) \, \mbox{ and } \, \theta_2=(z, \stackrel{g}{\ldots},z, sr, srz, z).$$Consequently, these actions are topologically equivalent. These surfaces form the hyperelliptic equisymmetric stratum $\theta_h,$ with the hyperelliptic involution being represented by $z$. We then conclude that  \begin{equation}\label{lapiz}\mathscr{M}_g(C_2^2, s=(0; 2^{g+3}), \theta_h)\end{equation}is a non-normal subvariety of $\mathscr{M}_g$ and that $\mathscr{S}_g$ consists of non-normal points of it. 

\s

Let $S$ be a member of \eqref{lapiz}. As proved in the previous case, $S$ is represented by  $$y^2=\Pi_{i=1}^{g+1}(x-a_i)(x-\tfrac{1}{a_i})$$where $a_ia_j \neq 1$ for each $i,j \in \{1, \ldots, g+1\},$ and the group $H_1$ is generated by the $$\langle z(x,y)=(x,-y), s(x,y)=(1/x,y/x^{g+1}) \rangle.$$ 

Assume now that $S$ belongs to  $\mathscr{S}_g.$ Then $r$ induces an automorphism  $\tilde{r}$ of $R=S/\langle z \rangle$ of order $g/2$ that keeps the set of points \begin{equation}\label{rojo}a_1, \ldots, a_{g+1}, \tfrac{1}{a_1}, \ldots, \tfrac{1}{a_{g+1}}\end{equation}invariant. Since the signature of the action of $\langle r \rangle$ on $S$ is $(0; 2^4, g, g)$ one sees that under the action of $\langle \tilde{r} \rangle \cong C_{g/2}$ the points \eqref{rojo} split into four long orbits and two fixed points. Consequently, up to a M\"{o}bius transformation, it can assumed that$$\tilde{r}(\xi)=\omega_{\frac{g}{2}} \xi \, \mbox{ and therefore } a_i= \omega_{\frac{g}{2}}^{i-1} a_1, \, a_{\frac{g}{2}+i}=\omega_{\frac{g}{2}}^{i-1} a_{\frac{g}{2}+i} \mbox { for } i=\{1, \ldots, \tfrac{g}{2}\}.$$ In addition, as $0$ and $\infty$ are the fixed points of $\tilde{r}$ we see that $a_{g+1}=0$ and $a_{g+1}=\infty.$

Set $a:=a_1^{g/2}$ and $b:=a_{g/2+1}^{g/2}$  The regular covering map induced by $\langle \tilde{r} \rangle$  $$R \to R'=R/\langle \tilde{r} \rangle \cong S/\langle r \rangle\,\mbox{ is given by } \,  \xi \mapsto \xi^{\frac{g}{2}}$$and therefore  the branch values of the composed map $S \to R'$ are $0, \infty$ (marked with $g$ each) and $a, 1/a, b, 1/b$ (marked with $2$ each). 

\s

The fact that $\langle  r \rangle$ is a normal subgroup of $G$ implies that the automorphism $s$ of $S$ induces an involution $\hat{s}$ of $R'$ in such a way that $S/G \cong R'/\langle \hat{s}\rangle.$ Clearly $\hat{s}(\xi)=1/\xi$ and $\pm 1$ are its fixed points. It then follows that the regular covering map induced by $\hat{s}$ $$R' \to R'/\langle \hat{s} \rangle \cong S/G \, \mbox{ is given by } \, \xi \mapsto \xi + 1/\xi,$$and the branch values of the composed map $S \to S/G$ are $\infty$ marked with $g$, $\pm 2$ marked with 2, and $a+1/a, b+1/b$ marked with 2.
 
\s

Thus, each member of $\mathscr{S}_g$ is  represented by the  curve \begin{equation*}\label{pz}y^2=x(x^{\frac{g}{2}}-a)(x^{\frac{g}{2}}-\tfrac{1}{a})(x^{\frac{g}{2}}-b)(x^{\frac{g}{2}}-\tfrac{1}{b})\end{equation*}for some distinct complex numbers $a,b \neq 0, \pm 1$ such that $ab \neq 1.$ In this model, the group of automorphisms  isomorphic to $G$ is generated by the transformations $$s(x,y)=(1/x, y/x^{g+1}) \mbox{ and } r(x,y)=(\omega_{g/2} x, \nu y)$$where $\nu$ satisfies $\nu^2=\omega_{g/2}$ and $\nu^{g/2}=-1.$

\s

Finally, if we take $a=-b=\eta$ where $\eta$ is a primitive fourth root of $-1$, then the equation above turns into 
$$y^2=x(x^{2g}+1),$$and therefore the Wiman curve of type II belongs to $\mathscr{S}_g,$ as desired.

\section{Proof of Theorem \ref{conj2}}
Let $g \geqslant 3$ be an odd integer. The semidirect product $$G=C_{2(g-1)} \rtimes C_{2}^2=\langle a, b, c : a^{2(g-1)}, b^2, c^2, (bc)^2, baba^{2-g}, caca^{-g} \rangle$$acts on a one-dimensional family $\mathscr{U}_g$ of Riemann surfaces of genus $g$ with signature $(0; 2^3,4).$ The existence of the family is guaranteed by the SKE \begin{equation}\label{pau}\theta=(b,bc,ab,bca^{-g}).\end{equation}Up to possibly finitely many exceptions, the full automorphism group of the members of $\mathscr{U}_g$ is isomorphic to $G.$ The equalities $$bc(a)bc=a^{-1}, \,\, (ac)^{g-1}=a^{g-1} \,\, \mbox{ and } \,\, b(ac)b=(ac)^{-1}$$show that $G$ has two  subgroups$$H_1=\langle a, bc \rangle  \mbox{ and }  H_2=\langle ac,b \rangle \mbox{ that are isomorphic to } \mathbf{D}_{2(g-1)}.$$As these groups are normal in $G$ we see that they are non-conjugate.

\s

Now, we determine the signature of the action of $H_1$ and $H_2$ on each $S \in \mathscr{U}_g.$ Since the groups $H_1$ and $H_2$ are normal, it is enough to compute the intersections $$\langle h \rangle \cap H_i \, \mbox{ for each } h \mbox{ appearing in \eqref{pau}}.$$We denote by $\pi, \pi_1$ and $\pi_2$ the regular covering maps given by the action of $G, H_1$ and $H_2$ respectively, and by $q_1, q_2, q_3,q_4$  the ordered branch values of  $\pi$. 

\s

{\bf a.} The $\pi$-fiber of $q_1$ consists of $4(g-1)$ points with $G$-stabiliser conjugate to $\langle b \rangle.$ They  yield two branch values of $\pi_2$ and does not yield branch values of $\pi_1.$

\s

{\bf b.} The $\pi$-fiber of $q_2$ consists of $4(g-1)$ points with $G$-stabiliser conjugate to $\langle bc \rangle.$ They yield two branch values of $\pi_1$ and does not yield branch values of $\pi_2.$

\s

{\bf c.} The $\pi$-fiber of $q_4$ consists of $4(g-1)$ points with $G$-stabiliser conjugate to $\langle a^gbc \rangle.$ They yield two branch values of $\pi_1$ and $\pi_2$ each.

\s

{\bf d.} The $\pi$-fiber of $q_3$ consists of $2(g-1)$ points with $G$-stabiliser conjugate to $\langle ab \rangle.$ Since $$(ab)^2=a^{g-1} \in H_i \mbox{ for } i=1,2,$$these points yield one branch value of $\pi_1$ and $\pi_2$ each.

\s

Hence, the signature of the action of both $H_1$ and $H_2$ on $S\in \mathscr{U}_g$ is $(0;2^5)$. Now, the proof of the theorem follows from the following claim.

\s

{\it Claim.} There is a unique topological class of actions of $$\mathbf{D}_{2(g-1)}=\langle r,s : r^{2(g-1)}, s^2, (sr)^2\rangle$$ in genus $g$ with signature $(0; 2^5).$

\s

We denote by $\theta=(g_1, g_2, g_3, g_4, g_5)$ the SKE representing an action of $\mathbf{D}_{2(g-1)}$  with signature $(0; 2^5).$  It is straightforward to see that among the elements $g_i$ exactly one of them must equal $z:=r^{g-1}.$ We may assume $g_5=z.$ In addition, if we write $g_i=sr^{n_i}$ then $$n_2-n_1+n_4-n_3 \equiv (g-1) \mbox{ mod } 2(g-1)$$and therefore among the $n_i$ exactly two of them are odd; we may assume that they are $n_2$ and $n_4$. After suitable conjugation coupled with an automorphism of the form $r \mapsto r^m$, we then have that $\theta$ is equivalent to the SKE $$(s, sr, sr^{n_3}, sr^{n_4}, z).$$Since the product of such elements must be trivial, we obtain 
$$\theta \cong \theta_{n_4}:= (s,sr,sr^{n_4+g}, sr^{n_4}, z).$$

Now, if $\Phi_3 \in \mathscr{B}$ is the braid transformation $(x_1, x_2, x_3, x_4) \mapsto (x_1, x_2, x_4, x_4^{-1}x_3x_4)$ then  $$\Phi_3 \circ \Phi_3 \cdot \theta_{n_4}=\theta_{n_4+2}$$Thus, $\theta \cong \theta_1=(s, sr, sr^{g+1}, sr, r^{g-1})$ and the claim follows.

\section{Proof of Theorem \ref{conj7}}
Let $n \geqslant 3$ be an odd integer and set $m=\tfrac{n-1}{2}$. Consider the group $$G=C_{2} \times C_{2n}=\langle a, b  : a^{2}, b^{2n}, abab^{-1} \rangle.$$The  Riemann-Huwitz formula together with the SKE \begin{equation*}\label{tarde}\theta=(a,b^n,b^2,\stackrel{m}{\ldots}, b^2,ab)\end{equation*}show that $G$ acts on a $m$-dimensional family $\mathscr{R}_g$ of Riemann surfaces  of genus $(n-1)^2$ with signature $(0; 2^2,n, \stackrel{m}{\ldots}, n,2n).$ Notice that, up to possibly finitely many exceptions, the full automorphism group of the members of $\mathscr{R}_g$ is isomorphic to $G.$  Clearly, the subgroups $$H_1=\langle b \rangle \,\, \mbox{ and } \,\, H_2=\langle ab^2 \rangle$$are isomorphic but non-conjugate in $G.$ We denote by $\pi, \pi_1$ and $\pi_2$ the regular covering maps given by the action of $G, H_1$ and $H_2$ respectively, and by $q_1, q_2, q_3^1,\ldots, q_3^m, q_4$  the ordered branch values of $\pi$. 

\s

{\bf a.} The $\pi$-fiber of $q_1$ consists of $2n$ points with $G$-stabiliser $\langle a \rangle.$ They  yield two branch values of $\pi_2$ with $H_2$-stabiliser $\langle a \rangle$ and does not yield branch values of $\pi_1.$

\s

{\bf b.} The $\pi$-fiber of $q_2$ consists of $2n$ points with $G$-stabiliser $\langle b^2 \rangle.$ They yield two branch values of $\pi_1$ with $H_1$-stabiliser $\langle b^n \rangle$ and does not yield branch values of $\pi_2.$

\s

{\bf c.} The $\pi$-fiber of $q_4$ consists of $2$ points with $G$-stabiliser $\langle ab \rangle.$ They yield one branch value of $\pi_i$ with $H_i$-stabiliser $\langle b^2 \rangle$ for $i=1,2.$

\s

{\bf d.} For each $i \in \{1, \ldots, m\}$, the $\pi$-fiber of $q_3^i$ consists of $4$ points with $G$-stabiliser $\langle b^2 \rangle.$ They yield two branch value of $\pi_i$ with $H_i$-stabiliser $\langle b^2 \rangle$ for $i=1,2.$

\s

Hence, the signature of the action of $H_1$ and $H_2$ on $S\in \mathscr{R}_g$ is $(0;2^2, n, \stackrel{n}{\ldots}, n)$ and the generating vector representing such actions have the form $$(b^n,b^n, g_1^1, g_1^2, \ldots, g_m^1, g_m^2, g_4) \, \mbox{ and } \, (a,a, h_1^1, h_1^2, \ldots, h_m^1, h_m^2, h_4)$$respectively, where $$\langle g_4 \rangle = \langle h_4 \rangle = \langle g_i^l\rangle =\langle h_i^l\rangle =\langle b^2 \rangle \, \mbox{ for each } i=1, \ldots, m \mbox{ and }l=1,2.$$

\s

The fact that $G$ is abelian implies that for each $P \in \pi^{-1}(q_4)$, the rotation constant of $P$ for $ab$ is given by $$\mbox{rot}(P, ab)=\omega_{2n} \, \mbox{ and therefore } \mbox{rot}(P, (ab)^2=b^2)=\omega_{2n}^2=\omega_n.$$See, for instance, \cite[\S3.1]{B22}. Now, by a result of Harvey, a SKE representing the action of $H_i$ can be obtained from the rotation numbers above. More precisely, following \cite[Theorem 7]{H71}, we see that $g_4=h_4=b^2$. 

\s

In a very similar way we deduce that $g_i^l=h_i^l=b^2$ as well.

\s

Hence, the proof now follows after noticing that the induced SKEs $$(b^n, b^n, b^2, \ldots, b^2) \, \mbox{ and } \, (a,a,b^2, \ldots, b^2)$$ are topologically equivalent.

\section{Proof of Theorem \ref{CFG}}
\subsection*{Generalised Fermat curves} Let $k,n \geqslant 2$ be integers such that $(k-1)(n-1)>2$.
We recall that a pair $(S,H)$ is called a {\it generalised Fermat pair} of type $(k,n)$ if $S$ is a compact Riemann surface endowed with a group of automorphisms $H$ isomorphic to $C_{k}^{n}$ in such a way that the signature of the action of $H$ is $(0; k, \stackrel{n+1}{\ldots}, k).$ If we write $$H=\langle a_1, \ldots, a_n\rangle \cong C_k^n$$then the existence of a generalised Fermat pair is guaranteed by the SKE $$\Theta_{k,n} :=(a_1, \ldots, a_{n},a_{n+1}) \mbox{ where } a_{n+1}:=(a_1\cdots a_n)^{-1}.$$In such a case, we say that  $S$ (respectively, $H$) is a {\it generalised Fermat curve} (respectively, a {\it generalised Fermat group}) of type $(k,n)$. Observe that the generalised Fermat curves of type $(k,n)$ form a $(n-2)$-dimensional  family in $\mathscr{M}_g$, where $$g=1+\tfrac{1}{2}k^{n-1}((n-1)(k-1)-2).$$

After considering a suitable M\"obius transformation, we may  assume that the branch values of the associated regular covering map $S \to S/H \cong \mathbb{P}^1$ are
\begin{equation}\label{bset}\infty, 0, 1,\lambda_{1},\ldots,\lambda_{n-2},\end{equation}
where $\lambda_{1}, \ldots, \lambda_{n-2} \in {\mathbb C} - \{0,1\}$ are pairwise distinct. 

\s

We now consider the irreducible and smooth algebraic curve
$$X_{\lambda_{1},\ldots,\lambda_{n-2}} :  \left\{ \begin{array}{ccccccc}
x_{1}^{k}&+&x_{2}^{k}&+&x_{3}^{k}&=&0\\
\lambda_{1}x_{1}^{k}&+&x_{2}^{k}&+&x_{4}^{k}&=&0\\
\, &\vdots &\,&\vdots&\,&\vdots&\,\\
\lambda_{n-2}x_{1}^{k}&+&x_{2}^{k}&+&x_{n+1}^{k}&=&0
\end{array}\right \} \subset \mathbb{P}^n
$$together with the linear automorphisms of ${\mathbb P}^{n}$ given by
$$\hat{a}_{j} : [x_{1}:\ldots:x_{n+1}] \mapsto [x_{1}:\ldots:x_{j-1}: \omega_{k} x_{j}: x_{j+1}:\ldots:x_{n+1}]$$for each $j=1,\ldots,n+1,$ where  $\omega_{k}=\mbox{exp}(\tfrac{2 \pi i}{k})$. Observe that$$\hat{a}_{1} \cdots \hat{a}_{n+1}=1 \mbox{ and } C_{k}^{n} \cong \langle \hat{a}_{1}, \ldots, \hat{a}_{n} \rangle =H_{0} \leqslant \mbox{Aut}(X_{\lambda_1, \ldots, \lambda_{n-2}}) .$$

In \cite{GHL} it was proved that $(X_{\lambda_1,\ldots,\lambda_{n-2}},H_{0})$ is a generalised Fermat pair of type $(k,n)$  and that there exists an isomorphism $$\Xi : S \to X_{\lambda_1,\ldots,\lambda_{n-2}} \mbox{ such that } H_0=\Xi H \Xi^{-1}.$$We can then identify $a_j$ with $\hat{a}_j.$

In terms of Fuchsian groups, the generalised Fermat pairs can be understood  as follows. If $\Gamma$ is a Fuchsian group such that ${\mathbb H}/\Gamma \cong S/H$  then $S\cong{\mathbb H}/\Gamma'$ and $H \cong \Gamma/\Gamma'$, where $\Gamma'$ stands for its derived subgroup of $\Gamma.$  In addition, $H$ is the unique generalised Fermat group of type $(k,n)$ of $S$ (see \cite{GHL} and \cite{HKLP}).

\subsection*{Proof of Theorem \ref{CFG}}
Let $J$ be a non-empty subset of $\{1,\ldots,n+1\}$ and write$$H_J=\langle  \{a_j : j \in J\} \rangle \leqslant H_0$$Note that if $|J| \in \{n,n+1\}$ then $H_{J}=H_{0}$. In the remaining cases we have that $$H_{J} \cong C_{k} \times \stackrel{|J|}{\ldots} \times C_{k} = C_{k}^{|J|}$$

Assume that $J_1$ and $J_2$ are two non-empty subsets of $\{1, \ldots, n+1\}$ such that $$|J_{1}|=|J_{2}| \leqslant n-2.$$

We claim that the actions of $H_{J_1}$ and $H_{J_2}$ on $X_{\lambda_{1},\ldots,\lambda_{n-2}}$ are topologically equivalent. Indeed, set $q_1=\infty, q_2=0, q_{3}=1$ and $q_{3+j}=\lambda_{j}$ for $1 \leqslant j \leqslant n-2$. Let $f$ be an  orientation-preserving self-homeomorphism of the Riemann sphere such that it preserves the  set $\{q_{1},\ldots,q_{n+1}\}$ and satisfies  \begin{equation}\label{hola7}f(\{q_{j}: j \in J_{1}\})=\{q_{l}: l \in J_{2}\}.\end{equation} 

Due to the fact that $\pi:S \to S/H$ is a homology branched covering map, the homeomorphism $f$ lifts to an orientation-preserving homeomorphism $$\widehat{f}:X_{\lambda_{1},\ldots,\lambda_{n-2}} \to X_{\lambda_{1},\ldots,\lambda_{n-2}} \mbox{ such that } \hat{f} H_{0} \hat{f}^{-1}=H_{0}.$$ Now, the property \eqref{hola7} of $f$ ensures that $\hat{f} H_{J_1} \hat{f}^{-1}=H_{J_2}$ and hence the actions are topologically equivalent, as claimed. 

\s

Let $\mbox{M\"{o}b}(\mathbb{C})$ denote the group of M\"{o}bius transformations. Let $K$ be the $\mbox{M\"{o}b}(\mathbb{C})$-stabiliser of the set of branch values \eqref{bset} of $S \to S/H.$
Since $H$ is a normal subgroup and $\Gamma'$ is a characteristic subgroup of $\Gamma$, there is a natural short exact sequence of groups
$$1 \to H \to {\rm Aut}(S) \to K \to 1.$$Hence, if  $K$ is the trivial group then $n \geqslant 4$ and $\mbox{Aut}(S) = H$.

\s

We now assume that $n \geqslant 4$ and $k \geqslant 2$, and choose $\lambda_{1}, \ldots, \lambda_{n-2}$ as before but satisfying that property that the set
$\{q_{1},\ldots, q_{n+1}\}$ has trivial $\mbox{M\"{o}b}(\mathbb{C})$-stabiliser. Observe that this property rules out only a positive codimensional subvariety of the parameters space. Now, if $S=X_{\lambda_{1},\ldots,\lambda_{n-2}}$ and $H=H_{0}=\langle a_{1}, \ldots,a_{n}\rangle$ then the aforementioned facts show that

\s

\begin{enumerate} 
\item the full automorphism group of $S$ agrees with $H$ (and thus is abelian), and 

\s

\item if $H_{J_{1}} \neq H_{J_{2}}$ then they are not conjugate in $\mbox{Aut}(S)$ and consequently their actions are not analytically equivalent.
\end{enumerate}

In particular, the actions of $$H_1=\langle a_1, \ldots, a_{n-1} \rangle \cong C_k^{n-1} \mbox{ and } H_2=\langle a_2, \ldots, a_{n} \rangle \cong C_k^{n-1}$$ on $S$ are topologically but not analytically equivalent, as desired.

\s

After considering the fact that the only non-trivial elements of $H_0$ acting with fixed points are the powers of $a_{1},\ldots,a_{n+1}$ and that every fixed point of a power of $a_{j}$ is also a fixed point of $a_{j}$ (see \cite{GHL}), 
it is a straightforward task to verify that the signature of the action of each $H_j$ on $S$ is $(0;k,\stackrel{k(n-1)}{\ldots},k)$. The induced actions are represented by the action of $H_1$ given by the SKE $$\theta_{k,n}:=(a_1, \stackrel{k}{\ldots}, a_1, \ldots, a_{n-1}, \stackrel{k}{\ldots}, a_{n-1}).$$ In conclusion \begin{equation}\label{final}\mathscr{M}_g(C_k^{n-1}, s=(0;k,\stackrel{k(n-1)}{\ldots},k), \theta_{k,n})\end{equation}is a non-normal subvariety of $\mathscr{M}_g$ of dimension $k(n-1)-3$, and a set of non-normal points consists of the generalised Fermat curves of type $(k,n)$.

\begin{rema}\mbox{}\label{uu}
\begin{enumerate}
\item The condition $n \geqslant 4$ in the theorem is needed as any collection of at most four points in $\mathbb{P}^1$ has non-trivial $\mbox{M\"{o}b}(\mathbb{C})$-stabiliser. 

\s

\item The theorem above says that (with $k=2$ and $n=4$) the subvariety \begin{equation}\label{llaves}\mathscr{M}_5(C_2^3, s=(0;2^6), \theta_{2,4}) \subset \mathscr{M}_5\end{equation}is non-normal, and a set of non-normal points is formed by the generalised Fermat curves of type $(2,4)$ (also known as classical Humbert curves). By considering the rule $X \to X/\langle a_1a_2 \rangle$, it can be seen that \eqref{llaves} is in bijective correspondence with the subvariety$$\mathscr{M}_3(C_2^2, s=(0;2^{6}), \theta_h)$$mentioned in the introduction as the first example of a non-normal subvariety of type \eqref{eq} associated to non-cyclic covers of $\mathbb{P}^1.$
\s

\item Similarly as pointed out in Remark \ref{telefono}, $X_{\lambda_1, \ldots, \lambda_{n-2}}$ does represent a non-normal point of \eqref{final} for all $\lambda_{1}, \ldots, \lambda_{n-2} \in {\mathbb C} - \{0,1\}$ pairwise distinct, despite the fact that,  for suitably chosen values of $\lambda_j$, the automorphism group of $X_{\lambda_1, \ldots, \lambda_{n-2}}$ does not have two non-conjugate subgroups acting in the same topological way. \end{enumerate}
\end{rema}

\section*{\bf Acknowledgement}
The first, third and fourth authors were partially supported by ANID Fondecyt Regular Grants 1230001, 1220099, 1230708 and 1230034. The second author was partially supported through the Fulbright Scholars Program, Fulbright Chile, and a Frank and Roberta Furbush Scholarship from Grinnell College.  The second author would also like to thank the faculty, staff, and students at Universidad de La Frontera and Universidad de Chile for their wonderful hospitality. The authors are grateful to the referees for their valuable comments and suggestions.

\subsection*{\bf Conflict of Interest} The authors have no conflict of interest to declare that are relevant to this article

\end{document}